\documentclass[3p,times,12pt,nopreprintline]{elsarticle}  

\usepackage{amssymb}
\usepackage{amsmath,bbm,units}
\usepackage{amsthm}
\usepackage[super]{nth}
\usepackage{graphicx}
\usepackage{caption}
\usepackage{newtxtext}
\usepackage[varvw]{newtxmath}

\usepackage[font+=footnotesize]{subcaption}
\usepackage{color}
\usepackage{tabularx,units}

\newtheorem{theorem}{Theorem}[section]
\newtheorem{remark}[theorem]{Remark}
\newtheorem{proposition}[theorem]{Proposition}
\newtheorem{lemma}[theorem]{Lemma}
\newtheorem{example}[theorem]{Example}

\newcommand{\bR}{{\mathbb{R}}}
\newcommand{\bN}{{\mathbb{N}}}

\newcommand{\cK}{{\mathcal{K}}}
\newcommand{\cL}{{\mathcal{L}}}
\newcommand{\cH}{{\mathcal{H}}}

\newcommand{\cP}{{\mathcal{P}}}
\newcommand{\cS}{{\mathcal{S}}}
\newcommand{\cW}{{\mathcal{W}}}

\begin{document}

\begin{frontmatter}

\title{Approximate Duals of B-splines for the Exact Representation of Splines on Coarse Knot Vectors}

\author[tud]{J.~St\"ockler\corref{cor}}
\ead{stoeckler@math.tu-dortmund.de}

\date{Version of \today}

\cortext[cor]{Corresponding author}

\address[tud]{Fakult\"at f\"ur Mathematik, Technische Universit\"at Dortmund, Vogelpothsweg 87, 44221 Dortmund, Germany}

\begin{abstract}
Approximate duals of B-splines were first used in \cite{Chuietal2004} for the purpose of constructing tight wavelet frames
on bounded intervals. They are splines with local support, whose inner product with a polynomial
in the spline space
provides the exact coefficient in the representation of the same polynomial in the B-spline basis.
This implies that the
 associated integral operator is a quasi-projection in the sense of \cite{Jia2004}.
Moreover, the approximation of smooth functions in Sobolev spaces
by this quasi-projection yields the optimal approximation order.
More recently, for applications in isogeometric analysis, the optimal approximation order should also be obtained for
functions in a slightly larger space, the so-called bent Sobolev space defined in \cite{Bazilevsetal2006,daVeigaetal2014}.
This requires the construction of enhanced approximate duals, whose corresponding integral operator
provides the exact representation of all spline functions with respect to a coarse knot vector,
using only few interior knots of the given knot vector.
Our analysis provides an explicit construction and hints for an efficient computation of enhanced approximate duals.
Explicit representations as linear combinations of B-splines are provided for $m=2$ and $m=3$, and two numerical examples
for splines of order $3\le m\le 6$ demonstrate the optimal approximation order for functions in bent Sobolev spaces.
\end{abstract}


\end{frontmatter}



\section{Introduction}
\label{sec:intro}

\noindent
Approximate duals of B-splines of order $m\in \bN$
were introduced in \cite{Chuietal2004} for the purpose of constructing tight wavelet frames
on bounded intervals. They received new attention in \cite{Dornischetal2016}
as substitutes for dual B-splines in the mortar finite-element method
and in connection with isogeometric analysis. Their main advantage in comparison with dual B-splines is their local support.
The sacrifice, again in comparison with dual B-splines, is the loss of an exact representation of all spline functions
by the associated integral operator. For a more precise explanation let $\cS_m(\Theta)$ be the space of spline functions of
order $m$ (degree $m-1$) with respect to the knot vector $\Theta\subset [a,b]$. Throughout this article,
$\Theta$ will denote an open knot vector; i.e. it has $m$-fold knots
at both endpoints $a,b$ and simple or multiple interior knots. The B-splines
\[  N_{m,k}(x) = (\theta_{k+m}-\theta_{k})~[\theta_{k},\ldots,\theta_{k+m}\mid (\cdot-x)_+^{m-1}],\quad 1\le k\le n,\]
define the standard basis $\Phi_m =(N_{m,1},\ldots,N_{m,n})$ of $\cS_m(\Theta)$.
Its dual basis
\[   \widetilde \Phi_m = (\widetilde N_{m,k};1\le k\le n) = \Phi_m \Gamma^{-1}, \]
where $\Gamma$ is the Gram matrix of $\Phi_m$, is useful to find the expression
\[  \widetilde \cK  f= \int_a^b f(y)\widetilde K(\cdot,y)\,dy,
\qquad \widetilde K(x,y)=\sum_{k=1}^n \widetilde N_{m,k}(y) N_{m,k}(x)
\]
for the orthogonal projection  $\widetilde K: L^2(a,b)\to\cS_m(\Theta)$.
Note that
the dual B-splines $\widetilde N_{m,k}$ have global support $[a,b]$.

As a substitute of dual B-splines,
  Chui et al. \cite{Chuietal2004} constructed {\em approximate duals}
\begin{equation}\label{eq:adual}
  \Phi_m \, S = ( N_{m,k}^{\rm ad};1\le k\le n),
\end{equation}
where $S$ is a symmetric positive-definite matrix of bandwidth $m-1$,  such that
\begin{equation} \label{eq:Ktheta}
   \cK  f= \int_a^b f(y) K(\cdot,y)\,dy,
\qquad  K(x,y)=\sum_{k=1}^n  N_{m,k}^{\rm ad}(y) N_{m,k}(x)
\end{equation}
reproduces all polynomials of degree $m-1$.
Hence, the reproducing property is retained for all polynomials in
the spline space, while the restriction to the local support of the splines $N_{m,k}^{\rm ad}$ sacrifices the reproducing property
for spline functions with ``active'' knots.
Linear operators $\cK:L^p(a,b)\to \cS_m(\Theta)$  with this property of polynomial reproduction
are called quasi-projections, see \cite{Jia2004}.
Since the reproduction of polynomials combined with local support are the fundamental properties for optimal
approximation rates in the Sobolev space $H^m(a,b)$, the integral operator $\cK$ in \eqref{eq:Ktheta} provides a
powerful tool for the approximation of smooth functions
by splines. The extension to higher dimension $d$ and rectangular domains $\hat \Omega\subset \bR^d$ can be done by
using the tensor product of univariate approximate duals.

More recently,  approximate duals were employed
 as Lagrange multipliers in domain decomposition methods in \cite{Dornischetal2016}.
 They replace the
dual B-splines used in   \cite{Belgacem1999,Wohlmuth2001} and lead to the implementation of weak continuity conditions
in sparse form.
Another application was presented in \cite{Dornischetal2024}, where approximate duals
(resp. the matrix $S$ in \eqref{eq:adual}) are used for mass lumping.
For similar applications in isogeometric analysis, the optimal approximation order of quasi-projections
should also be obtained for
functions in a slightly larger space, the so-called bent Sobolev space defined in \cite{Bazilevsetal2006,daVeigaetal2014}.
The bent Sobolev space of smoothness order $0\le s\le m$
associated with the knot vector $\Theta_0$ is
\begin{equation}\label{eq:bentSob}
    \cH^s([a,b];\Theta_0) = H^s(a,b)+ \cS_m(\Theta_0).
\end{equation}
Here, $\Theta_0$ is the knot vector which is used for the parameterization $X: [a,b]\to \bR^2$ of a curve $\gamma\subset \Omega$
in the physical domain $\Omega\subset\bR^2$. For example, for domain decomposition methods, $\gamma$ is
the interface between two subdomains of
$\Omega$.
The original definition in \cite{Bazilevsetal2006,daVeigaetal2014} is somewhat more technical, but the result in
\cite[Lemma 4.1]{daVeigaetal2014} shows that our definition in \eqref{eq:bentSob} gives the same space.

For functions $u\in H^s(\Omega)$, the pullback $\hat u_\gamma = u\circ X$ of the trace of $u$ on $\gamma$ is not always
in $H^{s-1/2}(a,b)$, due to the lack of high order smoothness of the parameterization $X$ of $\gamma$.
But it follows from the result in \cite[Lemma 4.1]{daVeigaetal2014}, that $\hat u_\gamma \in \cH^s([a,b];\Theta_0)$.
Therefore, in order to obtain optimal approximation rates for $\hat u_\gamma$ by a quasi-projection as in \eqref{eq:Ktheta},
the reproduction property should be valid for all functions
$f\in \cS_m(\Theta_0)$, not only for polynomials.

Throughout this article, we consider the typical situation where
$\cS_m(\Theta_0)$ is a subspace of $\cS_m(\Theta)$. We construct
{\em enhanced approximate duals} in $\cS_m(\Theta)$, whose corresponding integral operator
provides the exact representation of all splines in $\cS_m(\Theta_0)$. Typically,
$\Theta_0$ has only few interior knots of the given knot vector $\Theta$.
Our analysis provides an explicit construction and hints for an efficient computation of enhanced approximate duals.
Explicit represenations as linear combinations of B-splines are provided for $m=2$ and $m=3$,
and two numerical examples
for splines of order $3\le m\le 6$ demonstrate the optimal approximation order for functions in bent Sobolev spaces.

\section{Background on univariate B-splines, dual splines, and approximate dual splines}

\noindent
Our current work makes use of fundamental properties of univariate B-splines. Specifically, we use the recurrence relation for derivatives
and the explicit form for the representation of polynomials and truncated power functions in terms of the B-spline basis.
We
summarize the relevant properties as follows.

\begin{itemize}
\item We denote by $m\in\bN$ the order of the spline space $\cS_m(\Theta)\subset C[a,b]$, where
$\Theta=\{\theta_1,\ldots,\theta_{n+m}\}$ is an ``open'' knot vector. This means that $\Theta$ is a multiset with
$m$-fold knots at both endpoints, $\theta_1=\cdots=\theta_{m}=a$ and $\theta_{n+1}=\cdots=\theta_{n+m}=b$, and interior knots
$\theta_k$, $m+1\le k\le n$, of multiplicity at most $m-1$. The spline space has dimension $n$ and
a basis $\Phi_m = (N_{m,k};1\le k\le n)$ of
normalized
B-splines (of degree $m-1$), defined by
\[  N_{m,k}(x) = (\theta_{k+m}-\theta_{k})~[\theta_{k},\ldots,\theta_{k+m}\mid (\cdot-x)_+^{m-1}].\]
Instead of simple knot spans,  we use
\begin{equation}\label{eq:hmk}
  h_{j,k} = \frac{\theta_{k+j}-\theta_{k}}{j},\quad I_{j,k} = [\theta_{k},\theta_{k+j}],
\end{equation}
with obvious restrictions on the range of $j$ and $k$.

\item B-splines of higher order $m+\nu$ can be defined on the same knot vector $\Theta$, without increasing the multiplicity of the knots at the endpoints.  The corresponding spline space $\cS_{m+\nu}(\Theta)$ has dimension $n-\nu$
and basis $\Phi_{m+\nu}=(N_{m+\nu,k},~1\le k\le n-\nu)$.
Because the left endpoint $a$ is a knot of $N_{m+\nu,1}$
with multiplicity $m$, the function value and all derivatives up to order $\nu-1$ vanish at $x=a$. A similar statement holds at $x=b$ for
$N_{m+\nu,n-\nu}$.

\item The recursion for derivatives of B-splines can be written as
\begin{equation}\label{eq:Bsplineder}
   \frac{d}{dx} \Phi_{m+\nu+1} = \Phi_{m+\nu} D_{m+\nu}^T,\qquad \nu\ge 0,
\end{equation}
where the matrices $D_{m+\nu}$ of size $(n-\nu-1) \times (n-\nu)$ are given by
\begin{equation}\label{eq:Djmat}
 D_{m+\nu} = \left[ \begin{matrix} 1 &-1 \\ &1 &-1 \\ &&\ddots \\ & & &1 &-1 \end{matrix}\right]
~{\rm diag} \left( h_{m+\nu,1}^{-1},\ldots,h_{m+\nu,n-\nu}^{-1}\right).
\end{equation}

\item Marsden's identity for the spline space $\cS_m(\Theta)$ is
\begin{equation}\label{eq:Marsden1a} \nu!\binom{m-1}{\nu}(t-x)^{m-1-\nu} = \sum_{k=1}^n \phi_{m,k}^{(\nu)}(t) N_{m,k}(x),\qquad \phi_{m,k}(t) = \prod_{j=k+1}^{k+m-1} (t-\theta_j)
\end{equation}
for $0\le\nu\le m-1$.
E. Lee \cite[Corollary 2]{LeeETY1996} gives the formulation
\begin{equation}\label{eq:Marsden1b} \binom{m-1}{\nu} (t-x)^{m-1-\nu} = \sum_{k=1}^n \sigma_{m-1-\nu}(t-\theta_{k+1},\ldots,t-\theta_{k+m-1}) N_{m,k}(x)
\end{equation}
for $0\le \nu\le m-1$ with elementary symmetric polynomials
$$ \sigma_0\equiv 1, \quad \sigma_r (y_1,\ldots,y_{s}) = \sum_{1\le i_1<i_2<\ldots<i_r\le s} y_{i_1}\cdots y_{i_r} \quad\mbox{for}\quad s\ge r.$$

\item Similar formulas exist for truncated powers. If $\theta_\ell>\theta_{\ell-1}$ and $\theta_\ell$ appears in $\Theta$ with
multiplicity $\mu$, then by \cite[Corollary 1]{LeeETY1996}
\begin{equation}\label{eq:Marsden2a} \nu!\binom{m-1}{\nu}(\theta_\ell-x)_+^{m-1-\nu} = \sum_{k=1}^{\ell-m+\nu} \phi_{m,k}^{(\nu)}(\theta_\ell) N_{m,k}(x)\quad \mbox{for}\quad  0\le \nu\le \mu-1.
\end{equation}
A short argument leads to
\begin{equation}\label{eq:Marsden2b} \binom{m-1}{\nu} (\theta_\ell-x)_+^{m-1-\nu} =
\sum_{k=1}^{\ell-m+\nu} \sigma_{m-1-\nu}(\theta_\ell-\theta_{k+1},\ldots,\theta_\ell-\theta_{k+m-1}) N_{m,k}(x).
\end{equation}
The similarity with formulas \eqref{eq:Marsden1a} and \eqref{eq:Marsden1b} lies in the truncation of the power function
on the left-hand side and the truncation of the summation on the right-hand side.

{
For completeness of our exposition, we provide a short proof of \eqref{eq:Marsden2b}.
We let $0\le \nu\le \mu-1$ and take the $\nu$-th derivative in the first identity of the proof of Corollary 2 in \cite{LeeETY1996},
so
$$ \phi_{m,k}^{(\nu)}(s) = \sum_{j=0}^{m-1} \sigma_j(t-\theta_{k+1},\ldots,t-\theta_{k+m-1})~\frac{d^\nu}{ds^\nu} (s-t)^{m-1-j}.$$
This identity is true for all $t\in \bR$. When we choose $t=\theta_\ell$ and insert $s=\theta_\ell$, only the term with $j=m-1-\nu$ is nonzero.
Thus we obtain \eqref{eq:Marsden2b} from \eqref{eq:Marsden2a} and
$$ \phi_{m,k}^{(\nu)}(\theta_\ell) = \nu! \sigma_{m-1-\nu}(\theta_\ell-\theta_{k+1},\ldots,\theta_\ell-\theta_{k+m-1}).$$
}
\end{itemize}

The basis of dual B-splines $\widetilde \Phi_m = \Phi_m\Gamma^{-1}= (\widetilde N_{m,k};1\le k\le n)  $
is defined by means of the Gram matrix
\[  \Gamma= (\gamma_{j,k})_{1\le j,k\le n},\qquad \gamma_{j,k} = \int_a^b N_{m,j}(x)N_{m,k}(x)\,dx.
\]
The approximate duals of  Chui et al. \cite{Chuietal2004} have a similar definition as
\[
 \Phi_m^{\rm ad} = \Phi_m \, S = ( N_{m,k}^{\rm ad};1\le k\le n),
\]
where $S$ is a symmetric positive-definite matrix of bandwidth $m-1$, such that
\[
   {\rm supp}\,N_{m,k}^{\rm ad} =[\theta_{k_1},\theta_{k_2}],\quad k_1=\max(1,k-m+1),\quad k_2=\min(n+m,k+2m-1).
\]
The explicit formulation of $S$ is given in \cite[Section 5.2]{Chuietal2004} as
\begin{equation}\label{eq:defS}
  S =   U_0 + \sum_{\nu=1}^{m-1}    D_m\cdots D_{m+\nu-1} \,U_\nu\,D_{m+\nu-1}^T\cdots D_m^T
\end{equation}
with diagonal matrices $U_\nu={\rm diag}(u^{(\nu)}_k;1\le k \le n-\nu)$ and
\begin{equation}\label{eq:unu}
 u^{(\nu)}_k = h_{m+\nu,k}^{-1} F_\nu(\theta_{k+1},\ldots \theta_{k+m+\nu-1}),
\end{equation}
where $F_0 \equiv 1$ and
\[ F_\nu(x_1,\ldots,x_r)= (2^\nu \, \nu!)^{-1} \sum_{1\le i_1,\ldots,i_{\nu},j_1,\ldots,j_\nu\le r\atop {\rm distinct}}
   \prod_{\ell=1}^\nu (x_{i_{\ell}}-x_{j_\ell})^2\quad \text{for}\quad \nu\ge 1.
\]
Note that $F_\nu$ is a symmetric polynomial of degree $2\nu$ and the arguments of $F_\nu$ in \eqref{eq:unu} are
the interior knots of the B-spline $N_{m+\nu,k}$. Simplified expressions of $F_\nu$ in terms of centered moments of its
arguments $x_1,\ldots,x_r$ are also given in \cite{Chuietal2004}, together with Matlab code for the computation of
$u^{(\nu)}_k$.
Hence, $S$ is indeed symmetric and positive definite and has
bandwidth $m-1$, as stated in \eqref{eq:adual}.
 It was shown in \cite[Theorem 5.11]{Chuietal2004} that $S$ is the only
 symmetric matrix with bandwidth $m-1$, such that
 \[    K(x,y)=\Phi_m(y)~S ~\Phi_m(x)^T \]
 is the kernel of a quasi-projection. The bandwidth of $S$ implies that $K$ is localized
 around the diagonal; more precisely,
 if $x\in [\theta_k,\theta_{k+1})$ for some $m\le k\le n$, then
\begin{equation}\label{eq:kmsupport}
    K(x,y)=0\quad\text{if}\quad y\le \theta_{\max(1,k-2m+2)} \quad  \text{or}\quad y\ge \theta_{\min(n+m,k+2m-1)}.
\end{equation}
The definition of $S$ in \eqref{eq:defS} and the recurrence relation \eqref{eq:Bsplineder} give
\begin{equation}\label{eq:khq1}
   K(x,y) = \sum_{\nu=0}^{m-1} \sum_{k=1}^{n-\nu} u^{(\nu)}_k  \frac{\partial^{2\nu}}{\partial x^\nu \partial y^\nu}
   N_{m+\nu,k}(x)N_{m+\nu,k}(y).
\end{equation}

\section{Enhanced approximate duals and their quasi-projection}

\noindent
The quasi-projection onto $\cS_m(\Theta)$, which is defined by the kernel $K$ in \eqref{eq:khq1},
is useful for the approximation of smooth functions $f\in W^{s,p}(a,b)$.
Here we use the standard
notion of Sobolev spaces with a real smoothness parameter $s\ge 0$ and $1\le p\le \infty$ as defined
in \cite{Triebel1995}. The local support in \eqref{eq:kmsupport} and the polynomial reproduction provide
optimal error estimates of the form
\begin{equation}\label{eq:apporder}
  \left\|f-\int_a^b f(y)K(\cdot,y)\,dy\right\|_p \le  C h^s \|f\|_{s,p}
\end{equation}
for all functions $f\in W^{s,p}(a,b)$ and $0\le s\le m$, where $h$ on the right-hand side denotes a global
measure for
the local nature of the kernel $K$; a simple way to define it is
\[ h = \max\{ |x-y| : K(x,y)\ne 0\}.
\]

Applications of splines as finite elements in isogeometric analysis \cite{Bazilevsetal2006} require the approximation of
functions in slightly larger spaces, the bent Sobolev spaces defined in \eqref{eq:bentSob}, see
  \cite{Bazilevsetal2006,daVeigaetal2014}.
We provide a slightly more general form of these spaces for any $1\le p\le \infty$ by
 using the perspective based on the knowledge of the result in  \cite[Lemma 4.1]{daVeigaetal2014}.
The open knot vector $\Theta$ in $[a,b]$ is
 understood as the $h$-refinement of a knot vector $\Theta_0$, with $m$-fold knots at both endpoints of $[a,b]$ and
 a subset of the interior knots of $\Theta$.
For clarity when dealing with multiple knots, we define $\Theta_0$ by selecting $r$
indices
\[  m+1\le \ell_1 < \cdots <\ell_r \le n \]
of distinct interior knots of $\Theta$ such that $\theta_{\ell_j} > \theta_{\ell_j-1}$; i.e. $\theta_{\ell_j}$ is
the first occurence of this knot in $\Theta$. Then we choose multiplicities $\mu_j$
 less than or equal to the multiplicity of $\theta_{\ell_j}$ in $\Theta$ and obtain
\begin{equation}\label{eq:theta0}
   \Theta_0 = \{\underbrace{a, \ldots,a}_{m\text{-times}},\underbrace{\theta_{\ell_1}, \ldots,\theta_{\ell_1}}_{\mu_1\text{-times}},\ldots,\underbrace{\theta_{\ell_r}, \ldots,\theta_{\ell_r}}_{\mu_r\text{-times}},
   \underbrace{b, \ldots,b}_{m\text{-times}}\}.
\end{equation}
The bent Sobolev space of smoothness order $0\le s\le m$ associated with the knot vector $\Theta_0$ is defined by
\begin{equation}\label{eq:bentSob2}
    \cW^{s,p}([a,b];\Theta_0) = W^{s,p}(a,b)+ \cS_m(\Theta_0).
\end{equation}
The special case $p=2$ and $s\in \bN$ was considered in \cite{daVeigaetal2014}. Note that the sum is no direct sum, in particular
all polynomials of degree $m-1$ are elements of both spaces.

\begin{remark} {\rm In isogeometric analysis, splines (or NURBS) on tensor-product grids are used for the
 parameterization of the physical domain $\Omega$. Moreover, the finite element discretization for variational problems
 uses spaces of splines on refined grids ($h$-refinement). As a key observation, the pullback of smooth functions
 $u\in W^{s,p}(\Omega)$ to the parameter domain does not retain the same order of smoothness, due to the lack of high
 order smoothness of the parameterization of $\Omega$. Instead, the pullback is the element of a multivariate bent Sobolev space,
 defined by the sum of the standard Sobolev space on the parameter domain and the spline space for the parameterization.}
\end{remark}

In this article we provide a first step in order to extend the result for the approximation order in \eqref{eq:apporder}
to functions $f$ in the bent Sobolev space. For this purpose, we construct a new kernel $L$, as an enhancement of
the kernel $K$ for the quasi-projection, which reproduces all elements of $\cS_m(\Theta_0)$, not only polynomials.
More specifically, we will
define  $L\in \cS_m(\Theta)\otimes \cS_m(\Theta)$, such that
\[   \cL(p)= \int_a^b p(y) L(\cdot,y)\,dy  = p\quad  \text{for all}\quad p\in \cS_{m}(\Theta_0).
\]
The new kernel has the form
\begin{eqnarray}
  L(x,y) &=& K(x,y)+\sum_{j,k=1}^{n-m} u^{(m)}_{j,k}  \frac{\partial^{2m}}{\partial x^m \partial y^m}
   N_{2m,j}(x)N_{2m,k}(y)  \label{eq:Lm}  \\
   &=&  \Phi_{m}(y)~ S ~\Phi_{m}(x)^T+\frac{\partial^{2m}}{\partial x^m \partial y^m}
   \Phi_{2m}(y)~ U_{m} ~\Phi_{2m}(x)^T , \notag
\end{eqnarray}
which involves an additional matrix $U_m\in \bR^{(n-m)\times (n-m)}$ as compared to the kernel $K$.
In comparison with the definition of $S$ in \eqref{eq:defS}, we also let
\begin{equation}\label{eq:Sm1}
  S_L = S +  D_m^T\,\cdots \, D_{2m-1}^T\, U_{m}\,D_{2m-1}\,\cdots\,D_m
  \end{equation}
  and thus have
$$ L(x,y)= \Phi_{m}(y)~ S_L ~\Phi_{m}(x)^T.$$

Our goal is the definition of the matrix $U_m$ in the last summand of \eqref{eq:Sm1}, such that $\cL$ reproduces $\cS_m(\Theta_0)$.
In a first step, we characterize the reproduction of an arbitrary element $p= \Phi_m c\in \cS_m(\Theta)$
with coefficient vector $c\in \bR^n$. Note that we let $p$ be an element of the spline space associated with the
full knot vector $\Theta$ here.

\begin{lemma}\label{lem:reprodm}
Let $\Gamma$ denote the Gramian matrix of the B-spline basis $\Phi_m$ of $\cS_m(\Theta)$.
The reproduction of the spline $p= \Phi_m c\in \cS_m(\Theta)$  by the kernel $L$ in \eqref{eq:Lm} is equivalent to the linear equation
\begin{equation} \label{eq:reprodgeneralm}
  c^T (I_n-\Gamma S_L) = \vec 0.
\end{equation}
\end{lemma}

\begin{proof}
Inserting $p=\Phi_m c=c^T\Phi_m^T$  gives
\begin{eqnarray*}
 p(x)-\int_a^b p(y) L(x,y)\,dy &=& c^T \Phi_m(x)^T -\int _a^b c^T \,\Phi_m(y)^T\, \Phi_m(y) S_L ~\Phi_m(x)^T\,dy \\
 &=& c^T \left( I_n-\Gamma S_L\right) \Phi_m(x)^T.
 \end{eqnarray*}
 The linear independence of the B-splines in $\Phi_m$ gives the result.
\end{proof}

Because we wish to determine the summand of $S_L$ which contains the new matrix $U_m$, we rewrite
the condition \eqref{eq:reprodgeneralm} as
\begin{equation} \label{eq:reprodgeneralm2}
  c^T (I_n-\Gamma S) =  c^T \, \Gamma \, D_m^T\,\cdots \, D_{2m-1}^T\, U_{m}\,D_{2m-1}\,\cdots\,D_m.
\end{equation}

In order to achieve the reproduction of all splines in $\cS_m(\Theta_0)$,
we develop explicit formulas for the reproduction of all
 truncated powers of the form $(\theta_{\ell_j}-x)_+^{m-1-\nu}$ for $0\le \nu < \mu_j$, where $\ell_j$ and $\mu_j$ are
 the parameters in \eqref{eq:theta0}.
Let us first consider the reproduction of a single truncated power
$$p_{\ell,\nu}(x) := (\theta_\ell-x)_+^{m-1-\nu} = \Phi_m(x)c_{\ell,\nu}, \quad 0\le \nu\le \mu-1,$$
where $m+1\le \ell\le n$
is the index of an interior knot $\theta_\ell > \theta_{\ell-1}$ and $\mu$ is its multiplicity.
The coefficient vector $c_{\ell,\nu}$ is specified in  \eqref{eq:Marsden2b},
namely
\begin{equation}\label{eq:cellm}
c_{\ell,\nu;k}= \begin{cases}  \frac{1}{\binom{m-1}{\nu}} \sigma_{m-1-\nu}(\theta_\ell-\theta_{k+1} ,\ldots,\theta_\ell-\theta_{k+m-1})&
     \quad 1\le k\le \ell-m+\nu, \\
0,&\quad \ell-m+1+\nu\le k\le n.
\end{cases}
\end{equation}
Here, $\sigma_{m-1-\nu}$ is the elementary symmetric polynomial of degree $m-1-\nu$.
The row vector on the right-hand side of \eqref{eq:reprodgeneralm2} is simplified as follows.

\begin{lemma}\label{lem:leftsidem}
Assume  $\theta_\ell$ is an interior knot with multiplicity $1\le \mu\le m$  and $\theta_\ell> \theta_{\ell-1}$.
Let $0\le\nu\le \mu-1$ and
 $c_{\ell,\nu}$ as in \eqref{eq:cellm}. Then
$$ c_{\ell,\nu}^T \, \Gamma \, D_m^T\,\cdots \, D_{2m-1}^T = (m-1-\nu)!~\left( N_{2m,k}^{(\nu)}(\theta_{\ell}); ~1\le k\le n-m\right).$$
Here, $N_{2m,k}^{(\nu)}$ is the $\nu$-th derivative of the B-spline $N_{2m,k}$ of order $2m$ with respect to the given knot vector $\Theta$.
\end{lemma}

\begin{proof} The recursion for derivatives of B-splines can be written in terms of row vectors
$$   \left( N_{2m,k}^{(m)}; ~1\le k\le n-m\right) = \Phi_m D_m^T\,\cdots \, D_{2m-1}^T .$$
This gives
\begin{eqnarray*}
 \Gamma\,D_m^T\,\cdots \, D_{2m-1}^T &=& \int_a^b \Phi_m(y)^T \,\Phi_m(y)\, D_m^T\,\cdots \, D_{2m-1}^T\,dy \\[8pt]
   &=&
\int_a^b \Phi_m(y)^T \left( N_{2m,k}^{(m)}(y); ~1\le k\le n-m\right)\, dy.
\end{eqnarray*}
Multiplication from the left by the coefficient vector $c_{\ell,\nu}^T$ of the truncated power gives
$$  c_{\ell,\nu}^T \, \Gamma\,D_m^T\,\cdots \, D_{2m-1}^T = \int_a^b (\theta_{\ell}-y)_+^{m-1-\nu} \left( N_{2m,k}^{(m)}(y); ~1\le k\le n-m\right)\, dy .$$
Note that the order of all knots of the B-splines $N_{2m,k}$ cannot exceed $m$,
because $\Theta$ was defined as a knot vector for splines of order $m$. Therefore, the B-spline $N_{2m,k}$ and its first $m-1$ derivatives vanish at the endpoint $a$. We obtain by Taylor's formula
$$   \int_a^b (\theta_{\ell}-y)_+^{m-1-\nu}  N_{2m,k}^{(m)}(y)\, dy  = (m-1-\nu)!~ N_{2m,k}^{(\nu)}(\theta_{\ell}),\quad 1\le k\le n-m.$$
This completes the proof of the lemma.
\end{proof}

Another simplification of the identity \eqref{eq:reprodgeneralm2} is obtained by rewriting the row vector on the left-hand side
as follows.

\begin{lemma}\label{lem:rightsidem}
Assume  $\theta_\ell$ is an interior knot with multiplicity $1\le \mu\le m$
and $\theta_\ell> \theta_{\ell-1}$.
Let $0\le\nu\le \mu-1$ and
 $c_{\ell,\nu}$ as in \eqref{eq:cellm}. Then
\begin{equation}\label{eq:kappaellm}
 c_{\ell,\nu}^T  (I_n-\Gamma S) =  w_{\ell,\nu} \, D_{2m-1}\,\cdots \, D_{m},
 \end{equation}
 where the row vector $w_{\ell,\nu}\in\bR^{n-m}$ satisfies
 \begin{equation}\label{eq:kappaellm2}
    w_{\ell,\nu;k} = 0\quad\text{if}\quad k<\ell-2m+1+\mu\quad\text{or}\quad k>\ell-2.
 \end{equation}
\end{lemma}

\begin{proof} We let $v =  c_{\ell,\nu}^T  (I_n-\Gamma S) $ and obtain in a similar way as in the proof of Lemma \ref{lem:reprodm}
\begin{eqnarray*} q(x):=  v \, \Phi_m(x)^T &=&
 c_{\ell,\nu}^T   \Phi_m(x)^T - \int_a^b c_{\ell,\nu}^T\,\Phi_m(y)^T\,\Phi_m(y)\,S\, \Phi_m(x)^T\,dy  \\[8pt]
 &=& (\theta_\ell-x)_+^{m-1-\nu} -  \int_a^b (\theta_\ell-y)_+^{m-1-\nu} \,K(x,y)\,dy,
 \end{eqnarray*}
where $K$ is the kernel in \eqref{eq:khq1} for the reproduction of all polynomials in $\cS_m(\Theta)$.
For $x\in [a,\theta_{\ell-2m+1+\mu}]$, the property in \eqref{eq:kmsupport} implies that $K(x,y)=0$ for all
$$y\ge \theta_{\ell-2m+\mu+2m-1} = \theta_{\ell+\mu-1}=\theta_\ell.$$
Note that the last identity is valid due to the multiplicity of $\theta_\ell$.
Therefore, the truncation of $(\theta_\ell-y)_+^{m-1-\nu}$  at $\theta_\ell$ is not relevant for the integration.
Hence, the polynomial reproduction property gives
$$ q(x) = (\theta_\ell-x)^{m-1-\nu} -  \int_a^b (\theta_\ell-y)^{m-1-\nu} \,K(x,y)\,dy=0
$$
for all $x\in [a,\theta_{\ell-2m+1+\mu}]$.
Similarly, for $x\in [\theta_{\ell+2m-2},b]$, we have $K(x,y)=0$ for all $y\le \theta_\ell$ and
$$ q(x) = 0 -  0 =0.
$$
This proves that the spline $q = v \, \Phi_m^T$ has support in $[\theta_{\ell-2m+1+\mu},\theta_{\ell+2m-2}]$ and can be written as
$$ q= v \, \Phi_m^T = \sum_{j=\max(1,\ell-2m+1+\mu)}^{\min(n,\ell+m-2)} v_j N_{m,j}.$$
Consequently, all coefficients of the vector
$v =  c_\ell^T  (I_n-\Gamma S) $ with an index $j$, which lies outside the range of indices in the above summation, must be zero.

Next we show that $q$ has $m$ vanishing moments, that is
$$ \int_a^b x^\nu q(x)\,dx =\int_a^b x^\nu \left(p_{\ell,\nu}(x) -  \int_a^b p_{\ell,\nu}(y)\,K(x,y)\,dy\right)dx =0,\quad \nu=0,\ldots,m-1.$$
In fact, for every integrable function $f\in L_1(a,b)$, Fubini's theorem and the symmetry of $K$ give
$$ \int_a^b x^\nu \left( f(x)-\int_a^b f(y)K(x,y)\, dy\right)dx =
\int_a^b x^\nu f(x) \,dx -\int_a^b f(y)  \underbrace{\int_a^b x^\nu K(x,y)\, dx}_{\textstyle =y^\nu } \,dy =0.$$
Therefore, $q$ is the $m$-th derivative of a spline of order $2m$ with support in the same interval
$[\theta_{\ell-2m+1+\mu},\theta_{\ell+2m-2}]$. This means that
\begin{equation}\label{eq:kappaproof}
  q = v\, \Phi_m^T = \frac{d^m}{dx^m} ~\sum_{j=\max(1,\ell-2m+1+\mu)}^{\min(n-m,\ell-2)} w_j N_{2m,j},
\end{equation}
with a coefficient vector $w\in \bR^{n-m}$. The relation $v=  w\, D_{2m-1}\,\cdots \, D_{m}$ follows from
the recursion for derivatives of B-splines.
Note that the range of indices $j$ in \eqref{eq:kappaproof} corresponds to the property $w_k=0$ if
$ k<\ell-2m+1+\mu$ or $ k>\ell-2$ in \eqref{eq:kappaellm2}.
 \end{proof}

Next we proceed to the general case, where $\Theta_0\subset \Theta$ is defined as in \eqref{eq:theta0}.
Recall that $\theta_{\ell_j}$ is the first occurence of this knot in the knot vector $\Theta$.
Our goal is the simultaneous reproduction of
all truncated powers $(\theta_{\ell_j}-x)_+^{m-1-\nu}$ for $1\le j\le r$ and  $0\le \nu \le \mu_{j}-1$.
We let $\tilde r:= \sum_{j=1}^r \mu_{j}$ and obtain the following sufficient condition for the reproduction of these truncated powers.

\begin{proposition}\label{prop:findUm}
Let $A\in \bR^{\tilde r\times (n-m)}$ be the matrix with row vectors $a_{j,\nu}$ given by
\begin{equation}\label{eq:UmconditionA}
      a_{j,\nu} =(m-1-\nu)! \left(  N_{2m,k}^{(\nu)}(\theta_{\ell_j});~ 1\le k\le n-m\right),\quad 1\le j\le r,~~ 0\le \nu\le \mu_j-1,
\end{equation}
and let $B\in \bR^{\tilde r\times (n-m)}$ be the matrix with row vectors $w_{\ell_j,\nu}$ defined in \eqref{eq:kappaellm}.

Every solution $U_m$ of the equation
\begin{equation}\label{eq:AUmisB}
AU_m=B
\end{equation}
defines a kernel $L$ in \eqref{eq:Lm} which reproduces all polynomials of degree $m-1$ and all truncated powers
$(\theta_{\ell_j}-x)_+^{m-1-\nu}$ with $1\le j\le r$ and $0\le \nu\le \mu_j-1$.
\end{proposition}

\begin{proof} We first prove the reproduction of all
 truncated powers
\[ p_{\ell_j,\nu}(x)= (\theta_{\ell_j}-x)_+^{m-1-\nu} = \Phi_m(x) c_{\ell_j,\nu}
\]
by inserting the vectors  $c_{\ell_j,\nu}$ in identity \eqref{eq:reprodgeneralm2}. By Lemmas \ref{lem:leftsidem} and \ref{lem:rightsidem} we obtain the
equivalent identities
\[
  w_{\ell_j,\nu} \, D_{2m-1}\,\cdots \, D_{m}= (m-1-\nu)! \left(  N_{2m,k}^{(\nu)}(\theta_{\ell_j});~ 1\le k\le n-m\right)\, U_{m}\,D_{2m-1}\,\cdots\,D_m
\]
for all $1\le j\le r$ and $0\le \nu\le \mu_j-1$.
The definitions of $A$ and $B$ give the equivalent formulation
\[ B\, D_{2m-1}\,\cdots \, D_{m}= A\, U_{m}\,D_{2m-1}\,\cdots\,D_m.
\]
Therefore, the condition $AU_m=B$ is sufficient for the reproduction of all
 truncated powers as claimed in the proposition.

In addition, we show that the polynomial reproduction property of the kernel $   K(x,y) = \Phi_m(y) S \Phi_m(x)^T$
is not destroyed by adding the last summand in \eqref{eq:Sm1} to $S$.
 In other words, we show that
\begin{equation}\label{eq:polrepro}
  \int_a^b p(y) \Phi_m(y)~D_m^T\cdots D_{2m-1}^T\, U_m\,D_{2m-1}\cdots D_m \Phi_m(x)^T \,dy =0
\end{equation}
for all $p\in \cP_{m-1}$. The recursion for the derivatives of B-splines gives
$$ \Phi_m(y)~D_m^T\cdots D_{2m-1}^T=  \left( N_{2m,k}^{(m)}; ~1\le k\le n-m\right).$$
Recall that the multiplicity of each knot cannot exceed $m$. Therefore, the $m$-th derivative $ N_{2m,k}^{(m)}$ is well defined and
has $m$ vanishing moments.
This implies that the integral in \eqref{eq:polrepro} vanishes, which was to be shown.
\end{proof}

Finally, we can show that the matrix equation $AU_m=B$ in \eqref{eq:AUmisB} is always solvable, without any restriction on the
knot vector $\Theta$, the selection of interior knots $\theta_{\ell_1},\ldots,\theta_{\ell_r}$ and their multiplicities.
Note that $A$ has fewer rows than columns: the number $\tilde r$ of rows is less than or equal to the number $n-m$ of columns,
because $n-m$ is the sum of the multiplicities of all interior knots of $\Theta$.

\begin{lemma}\label{lem:rightinverse}
The matrix $A$ in \eqref{eq:AUmisB} has full rank $\tilde r$. Furthermore, every right-inverse $R$ of $A$
provides a solution $U_m=RB$ of \eqref{eq:AUmisB}.
\end{lemma}

\begin{proof} The second assertion is trivial, as we have $ARB = B$ for every right-inverse $R$ of $A$.
For the first assertion, we use results about total positivity of spline collocation matrices from \cite{deBoor1985}.
Note that the row vectors
$$ a_{j,\nu} =(m-1-\nu)! \left(  N_{2m,k}^{(\nu)}(\theta_{\ell_j});~ 1\le k\le n-m\right),\quad 1\le j\le r,~~ 0\le \nu\le \mu_j-1,$$
in \eqref{eq:UmconditionA} are rows of the B-spline collocation matrix
which is used for Hermite interpolation with B-splines $N_{2m,k}$, $1\le k\le n-m$. Here, we have a special situation where
 the nodes for interpolation are the knots $\{\theta_{\ell_j},1 \le j\le r\}\subset \Theta$
 and the multiplicities $\mu_j$
 of these nodes are less than or equal to the multiplicities in the knot vector $\Theta$.
 For clarity, we specify the set of interpolation nodes
 including their multiplicities by an ordered set
 \begin{equation}\label{eq:intnodes}
   X=\{x_1\le \cdots \le x_{\tilde r}\}\quad \text{with} \quad   x_i=\theta_{\ell_j}\quad\text{for}\quad
     \sum_{s=1}^{j-1}\mu_s+1 \le i \le \sum_{s=1}^{j}\mu_s,
 \end{equation}
 so that we can write the matrix $A$ as
 \begin{equation} \label{eq:Aexplicit}
  A =\left[(m-1-\nu_i)! ~  N_{2m,k}^{(\nu_i)}(x_i)\right]_{1\le i\le \tilde r;~1\le k\le n-m},\quad \text{where}\quad
 \nu_i=\#\{j: j<i~\text{with}~x_j=x_i\}.
 \end{equation}
The typical side conditions for Hermite interpolation with B-splines
 are satisfied. In particular, for each interpolation node $\theta_{\ell_j}$,
 \begin{enumerate}
 \item the sum of the multiplicities as an interpolation node and a spline knot
 does not exceed the order $2m$ of the B-splines,
 \item the rows $a_{j,\nu}$ associated with $\theta_{\ell_j}$
 involve consecutive derivatives  $N_{2m,k}^{(\nu)}$ with $0\le \nu\le \mu_j-1$.
 This explains that the condition (4.12) in \cite{deBoor1985} is satisfied.
 \end{enumerate}
 Based on these observations about the structure of the matrix $A$,
 we employ the following result of \cite[Theorem 3]{deBoor1985}: every square submatrix
 $$ A(k_1,\ldots,k_{\tilde r})\quad\text{consisting of the columns $k_1,\ldots,k_{\tilde r}$ of $A$}
$$
satisfies
 $$
 \det A(k_1,\ldots,k_{\tilde r}) \ge 0, $$
and strict positivity holds if and only if the interpolation nodes in \eqref{eq:intnodes} satisfy
\begin{equation}\label{eq:SWcondition}
    N_{2m,k_i}(x_i)>0,\qquad 1\le i\le \tilde r.
\end{equation}
The last condition \eqref{eq:SWcondition} is known as the Schoenberg-Whitney condition for spline interpolation.
A suitable selection of column indices $k_1,\ldots,k_{\tilde r}$ is described in the next result.
This completes the proof that $A$ has full rank, and thus there exist right inverses of $A$.
\end{proof}

\begin{proposition}\label{prop:Um2}
For each $1\le j\le r$ we define the set
\begin{equation}\label{eq:A0col1}
      K_j := \{  \ell_j-m ,\ldots,\ell_j-m+\mu_j-1 \} \subset \{1,\ldots,n-m\}
\end{equation}
of cardinality $\mu_j$. These sets are disjoint, and with
\begin{equation}\label{eq:A0col2}
 K:= \{k_1< \cdots< k_{\tilde r}\} = \bigcup_{j=1}^r K_j,
\end{equation}
the matrix $A_0:=A(k_1,\ldots,k_{\tilde r})$ satisfies $\det A_0>0$.
\end{proposition}

\begin{proof}
First, we note that $\ell_j-m\ge 1$ and $\ell_j-m+\mu_j-1\le n-m$, because $\theta_{\ell_j}=\theta_{\ell_j+\mu_j-1}$
is an interior knot of
$\Theta$. Moreover, the definition of $\Theta_0$ in \eqref{eq:theta0}
shows that $\theta_{\ell_{j+1}}>\theta_{\ell_j}=\theta_{\ell_j+\mu_j-1}$. Thus we have
$$  \min (K_{j+1}) -\max(K_j)=  \ell_{j+1}-m - (\ell_j-m+\mu_j-1) = \ell_{j+1} - (\ell_j +\mu_j - 1) >0.
$$
This shows that the sets $K_j$ do not overlap.
Their union defines a set $K$ of cardinality $\tilde r$, which is a subset of the column indices of $A$.
The corresponding submatrix $A_0$ is a square matrix. By comparing the notations in \eqref{eq:intnodes} and \eqref{eq:A0col2},
we see that
$$ x_i = \theta_{\ell_j}\quad\text{for all}\quad k_i\in K_j,$$
and this shows that
$$ N_{2m,k_i}(x_i)>0,\qquad 1\le i\le \tilde r.$$
Therefore, condition \eqref{eq:SWcondition} is satisfied and $\det A_0>0$ is proved.
\end{proof}

The next two examples describe possible choices of right inverses of the matrix $A$ in \eqref{eq:intnodes}
for the construction of the matrix $U_m$.

\begin{example}\label{example:MPinverse}
\rm
 The Moore-Penrose (MP) pseudo-inverse  of $A$ is a right-inverse. It is given by
$$ R=A^T(AA^T)^{-1} \in \bR^{(n-m)\times \tilde r}.$$
Note that $A$ is sparse, because at most $2m-1$ nonzero entries appear in each row.
Furthermore, in many practical situations, we have small values of $\tilde r$, say $\tilde r\le 5$,
in comparison with large values of $n$, so that the computation of $R$ is numerically inexpensive.
It is particularly simple, if the condition
\begin{equation}\label{eq:gapells}
 \ell_j\le \ell_{j+1}-2m +1 \quad\text{for all}\quad 1\le j\le r
\end{equation}
is satisfied, i.e. the selected interior knots have gaps of size at least $2m-1$.
Then the non-zero entries of $A$ appear in $r$ blocks, with size $\mu_j \times (2m-1)$ for $1\le j\le r$, whose row and column indices do not overlap. Therefore, the MP pseudo-inverse of $A$ is obtained by transposition of $A$
and blockwise pseudo-inverses. We include a numerical example with $m=2$, $\Theta=[0,0,1:7,8,8,9:11,12,12]$,
and $r=2$, choosing $\ell_1=5$ (simple knot $\theta_5=3$) and $\ell_2=10$ (double knot $\theta_{10}=8$). Then we compute
values of cubic B-splines at $x=3$ and $x=8$ and their first derivatives at $x=8$ to obtain
\[ A=\frac{1}{6} ~\left[\begin{array}{cccccccccccc}
                                0&1&4&1&0&0&0&0&0&0&0&0\\
                                0&0&0&0&0&0&0&3&3&0&0&0\\
                                0&0&0&0&0&0&0&-9&9&0&0&0
                                \end{array}\right].
\]
The MP pseudo-inverse is
\[  R=\frac{1}{3} ~\left[\begin{array}{cccccccccccc}
                                0&1&4&1&0&0&0&0&0&0&0&0\\
                                0&0&0&0&0&0&0&3&3&0&0&0\\
                                0&0&0&0&0&0&0&-1&1&0&0&0
                                \end{array}\right]^T.
\]
\end{example}

\begin{example}\label{example:A0}
\rm
Another type of right-inverse of $A$ is obtained as an application of the result in Proposition \ref{prop:Um2}.
We select $\tilde r$ columns of $A$ with indices $k_1,\ldots,k_{\tilde r}$ in \eqref{eq:A0col2}.
The corresponding matrix
\begin{equation}\label{eq:A0}
   A_0=A(k_1,\ldots,k_{\tilde r})
\end{equation}
is invertible. Then we define a right inverse $R$ of $A$
by putting the rows of $A_0^{-1}$ in positions $k_1,\ldots,k_{\tilde r}$ of $R$ and filling the remaining rows by zeros.
In the previous example, this gives column indices $k_1=3$, $k_2=8$, $k_3=9$ and an even sparser right inverse
\[  R= \frac{1}{6} ~\left[\begin{array}{cccccccccccc}
                                0&0&9&0&0&0&0&0&0&0&0&0\\
                                0&0&0&0&0&0&0&6&6&0&0&0\\
                                0&0&0&0&0&0&0&-2&2&0&0&0
                                \end{array}\right]^T.
\]
The effect of this procedure is most visible in the special case where
all selected interior knots $\theta_{\ell_j}$ are simple ($\mu_j=1)$ and
separated by gaps as in \eqref{eq:gapells}. Then the square submatrix $A_0$ in Proposition \ref{prop:Um2} is
a diagonal matrix with diagonal entries $N_{2m,\ell_j-m}(\theta_{\ell_j})>0$, $1\le j\le r$.
 Note that the choice of column indices of $A$
is such that the diagonal entries of $A_0$ are values of B-splines at their ``central'' knot.
In this case, inversion of $A_0$ is trivial and numerically stable.
\end{example}

\begin{remark} \rm
The following computational aspects are useful for an implementation to solve the matrix equation $AU_m=B$ in
\eqref{eq:AUmisB}.

\begin{itemize}
\item
The entries of $A$ (or $A_0$) can be efficiently computed using the Cox-de Boor algorithm, see \cite{deBoor2001}.
\item
We noticed in our numerical experiments up to $m=6$,
that the matrix $A_0$ has a small condition number $\le 10$ after division
of each row $a_{j,\nu}$ by its maximal entry in absolute value (diagonal preconditioning).
This preconditioning is especially effective if knots with higher multiplicities occur.

\item
The matrix $B$
is composed of row vectors $w_\ell$ in \eqref{eq:kappaellm}.
They are the result of stripping off the product of matrices $D_{2m-1}\cdots D_m$ from the vector
\begin{equation}\label{eq:vell}
 v_\ell = c_\ell^T(I_n-\Gamma S).
\end{equation}
Therefore, $w_\ell$  can be obtained by the right multiplication
\begin{equation}\label{eq:wfromv}
   w_\ell= v_\ell D_m^+\cdots D_{2m-1}^+, \quad D_j^+=  {\rm diag}\left(h_{j,k}\right)_{1\le k\le n+m-j} ~
\begin{pmatrix} 1& 1 & \cdots &1\\  0 & 1 &\cdots &1 \\ \vdots & \ddots & \ddots &\vdots \\
 0 &\cdots & 0 & 1 \\   0 &  \cdots& 0 &0 \end{pmatrix}.
\end{equation}
Note that right-multiplication of a row vector by an upper triangular matrix of ones is a cumulative sum of its entries, so no explicit
vector-matrix multiplication is required
(e.g. {\tt cumsum}-command in Matlab).

Finally, we do not need to compute all entries of the Gramian matrix $\Gamma$ in order to find $v_\ell$ in \eqref{eq:vell}.
Although $v_\ell$ has more nonzero entries than $w_\ell$, only its portion $\widetilde v_\ell:= (v_{\ell,k},~\ell-2m+2\le k\le \ell-2)$ is relevant
for the definition of the entries of $w_\ell$.
Moreover, all entries of the vector $c_\ell$ in \eqref{eq:cellm} with index $i\ge \ell-m+1$ are zero.
Since $\Gamma$ and $S$ have bandwidth $m-1$, only the entries of $\Gamma$ in positions $(i,j)$ with
$$  \max(1, \ell-4m+4) \le i \le \ell-m,\qquad  \max(\ell-3m+3,1) \le j\le \ell-1
$$
are relevant for the computation of $\widetilde v_\ell$.
\end{itemize}
\end{remark}

Section 5 contains explicit formulas for the matrix $B$ in Proposition \ref{prop:findUm} for the cases $m=2$ and $m=3$
with simple knots.

\section{Approximation of functions in bent Sobolev spaces}

\noindent
We explain by a numerical example that enhanced approximate duals provide the optimal approximation
order for functions in bent Sobolev spaces $\cH^m([a,b],\Theta_0) = H^m(a,b)+\cS_m(\Theta_0)$.
So we choose $p=2$ and $s=m$ in \eqref{eq:bentSob2}. A detailed error analysis
for the general case is postponed to a forthcoming article.

Our setup is the same as in Example 9.1.2 in \cite{DornischStoeckler2021}, where dual B-splines are used
  as Lagrange multipliers in a problem of linear elasticity. Let $\Omega = (0,1)^2$ and $u(x,y)=\sin(3x)\sin(2y)$.
 We split $\Omega$ into two subdomains by the quadratic spline curve $\gamma$ with parameterization
 \[ X:(0,1)\to \Omega,\qquad X(t) = N_{3,1}(t){0.5\choose 0} +	N_{3,2}(t)	{0.6\choose 0.3}
    + 	N_{3,3}(t)	{0.4\choose 0.7} +  N_{3,4}(t)	{0.5\choose 1}
    \]
of order  $3$ with
knot sequence ${\Theta_0}=\{0,0,0,0.5,1,1,1\}$. 
The point $X(0.5)=(0.5,0.5)$
corresponds to the parameter value $t=0.5$, where the curve has only $C^1$-smoothness. When we represent the same curve
by splines of higher order $m\ge 3$, the corresponding knot vector must be adapted to
\[ \Theta_{0,m} = \{\underbrace{0,\cdots,0}_{m},\underbrace{0.5,\cdots,0.5}_{m-2},\underbrace{1,\cdots,1}_{m}\} \]
and the new control points are obtained by simple convex combinations, described as degree elevation in \cite[Chapter 8]{Farin2002}.

Whereas the function $u$ is an element of $H^s(\Omega)$ for every $s\ge 0$, its pullback $\hat u=u\circ X$ to the interface $\gamma$
is an element
of $H^s(0,1)$ only for $0\le s<5/2$, due to the reduced smoothness of $X$ at $t=0.5$. On the other hand, when $m\ge 3$, then
$\hat u$ is an element
of the bent Sobolev space
\[ \cH^m((0,1);\Theta_{0,m}) = H^m(0,1)+ \cS_m(\Theta_{0,m}). \]
This can be seen as follows. The multiplicity of the knot $0.5$ in $\Theta_{0,m}$ is $m-2$, so
the truncated power $(\cdot-0.5)_+^2$ is an element of $\cS_m(\Theta_{0,m})$. Therefore, the
 jump in the second derivative of $\hat u$ at $t=0.5$
can be absorbed by the spline
\[ s_1=   c_2(\cdot-0.5)_+^2 \in \cS_m(\Theta_{0,m}),\quad c_2=\frac{1}{2}(\hat u^{(2)}(0.5+)-\hat u^{(2)}(0.5-)) ,\]
so that $\hat u-s_1\in C^2[0,1]$. Continuing in the same way, we can define the spline
\[  s = \sum_{\nu=2}^{m-1}   c_\nu (\cdot-0.5)_+^\nu\in \cS_m(\Theta_{0,m}),\]
such that $f:=\hat u-s\in C^{m-1}[0,1]$. Since $f$ is piecewise $C^\infty$, we easily obtain $f\in H^m(0,1)$.
The above construction of the spline $s\in \cS_m(\Theta_{0,m})$ is the essence of the proof of Lemma 4.1 in \cite{daVeigaetal2014}.

\begin{figure}[ht]
	\centering
	\includegraphics[width=0.4\textwidth]{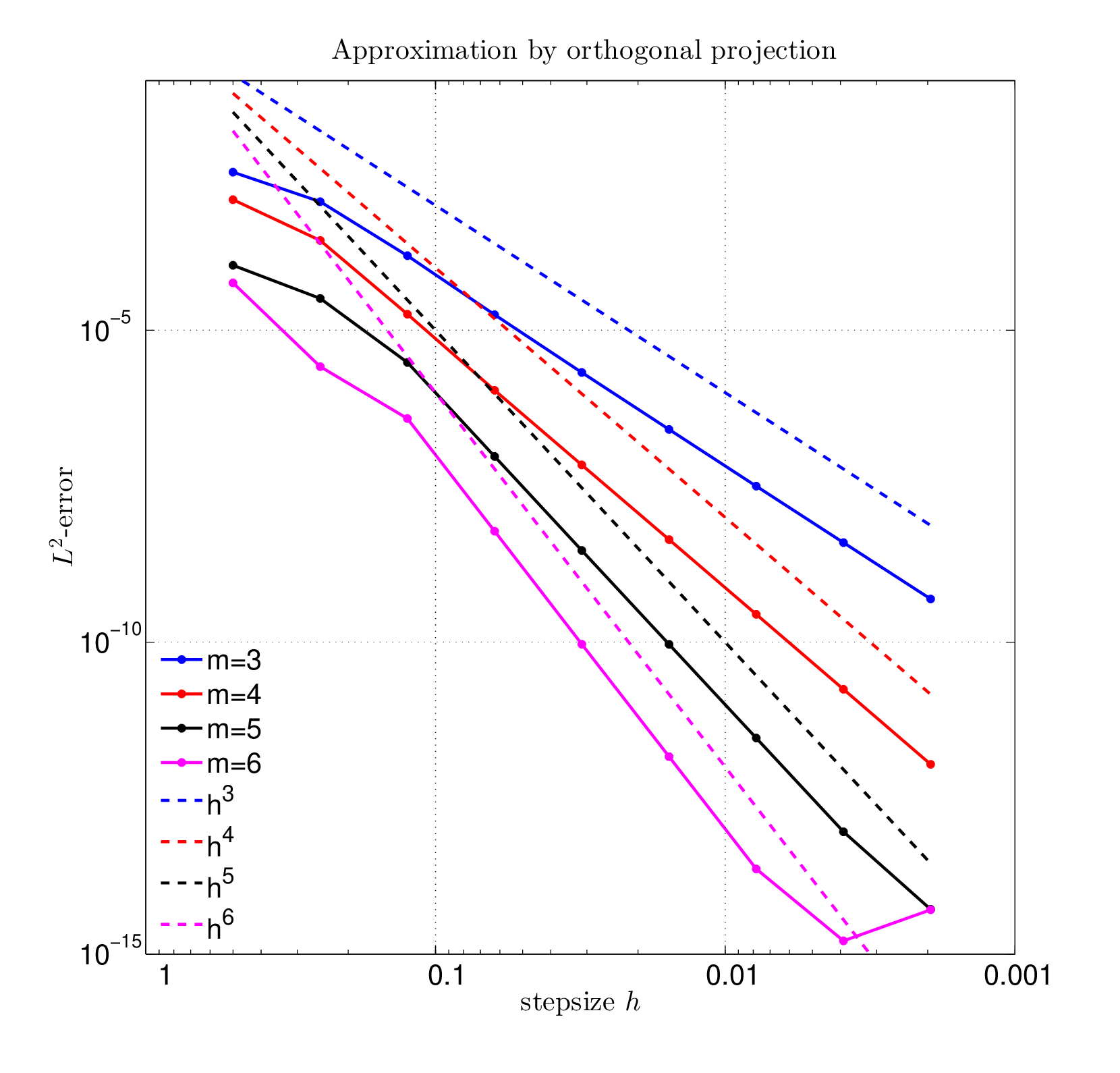}\quad
    \includegraphics[width=0.4\textwidth]{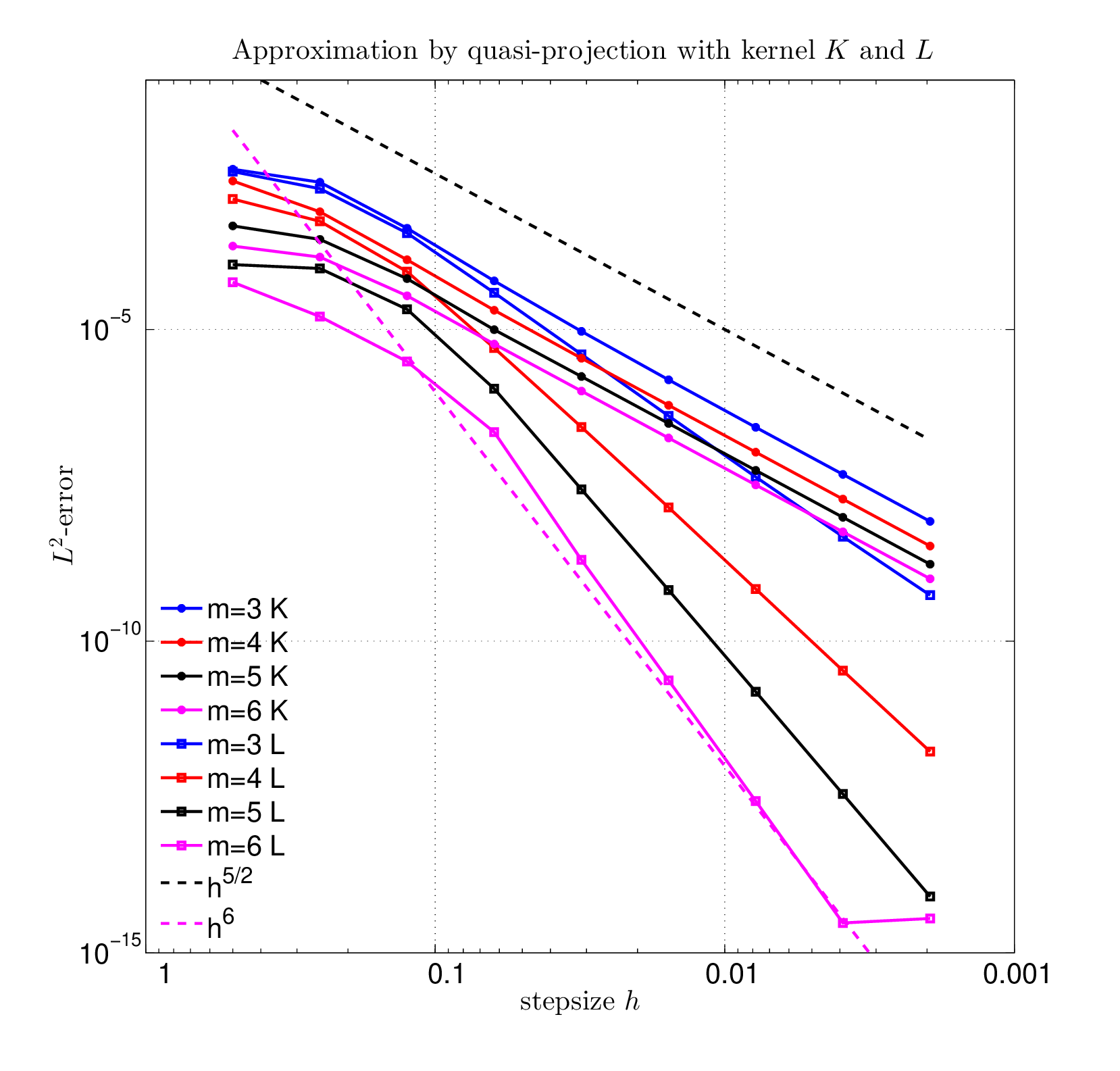}
  \caption{$L^2$-error for approximation of $\hat u$ by splines of order $3\le m\le 6$, with orthogonal projection (left)
  and quasi-projections by means of the kernels $K$ and $L$ (right).
  \label{fig:approxu}}
\end{figure}

Our numerical results are obtained using matlab.
First, we consider the approximation of $\hat u$ by splines in $\cS_m(\Theta_{h,m})$ of order $3\le m\le 6$.
The knot vector $\Theta_{h,m}$ is defined as a refinement of $\Theta_{0,m}$ by inserting equidistant knots with stepsize $h=\tfrac{1}{2N}$
into both intervals $[0,0.5]$ and $[0.5,1]$. These new knots have multiplicity $1$.
We show in Figure \ref{fig:approxu} that the quasi-projection of $\hat u$ by means of the kernel $L$ in \eqref{eq:Lm} provides the
same approximation order $O(h^m)$ for the $L^2$-error as the orthogonal projection into this spline space. This supports our claim that
the enhanced
approximate duals provide a good substitute of dual B-splines, have local support and obtain the optimal approximation order
for functions in the bent Sobolev space.
On the other hand, the quasi-projection of $\hat u$ by means of the kernel $K$ in \eqref{eq:khq1} only reaches
the approximation order $O(h^{2.5})$, which corresponds to the sharp upper bound of the Sobolev regularity of $\hat u$,
even for higher order $m=6$.

Next we extend our discussion to the aforementioned application in isogeometric analysis, where dual B-splines
are used as Lagrange multipliers on interface curves
in domain decomposition methods. Here the approximation of the normal derivative $g:=\frac{\partial u}{\partial n}$
along $\gamma$ should also provide the optimal approximation order. The observation in \cite[Example 9.1.2]{DornischStoeckler2021}
made clear that the pullback $\hat g$ of this function is in $H^s(0,1)$ only for $0\le s<3/2$,
but in the bent Sobolev space
\[ \cH^m((0,1);\Theta_{0,m}^*) = H^m(0,1)+ \cS_m(\Theta_{0,m}^*)\quad \text{for}\quad m\ge 3, \]
where $\Theta_{0,m}^*$ has an increased multiplicity of the interior knot $0.5$, so
\[ \Theta_{0,m}^* = \{\underbrace{0,\cdots,0}_{m},\underbrace{0.5,\cdots,0.5}_{m-1},\underbrace{1,\cdots,1}_{m}\}. \]

\begin{figure}[ht]
	\centering
	\includegraphics[width=0.4\textwidth]{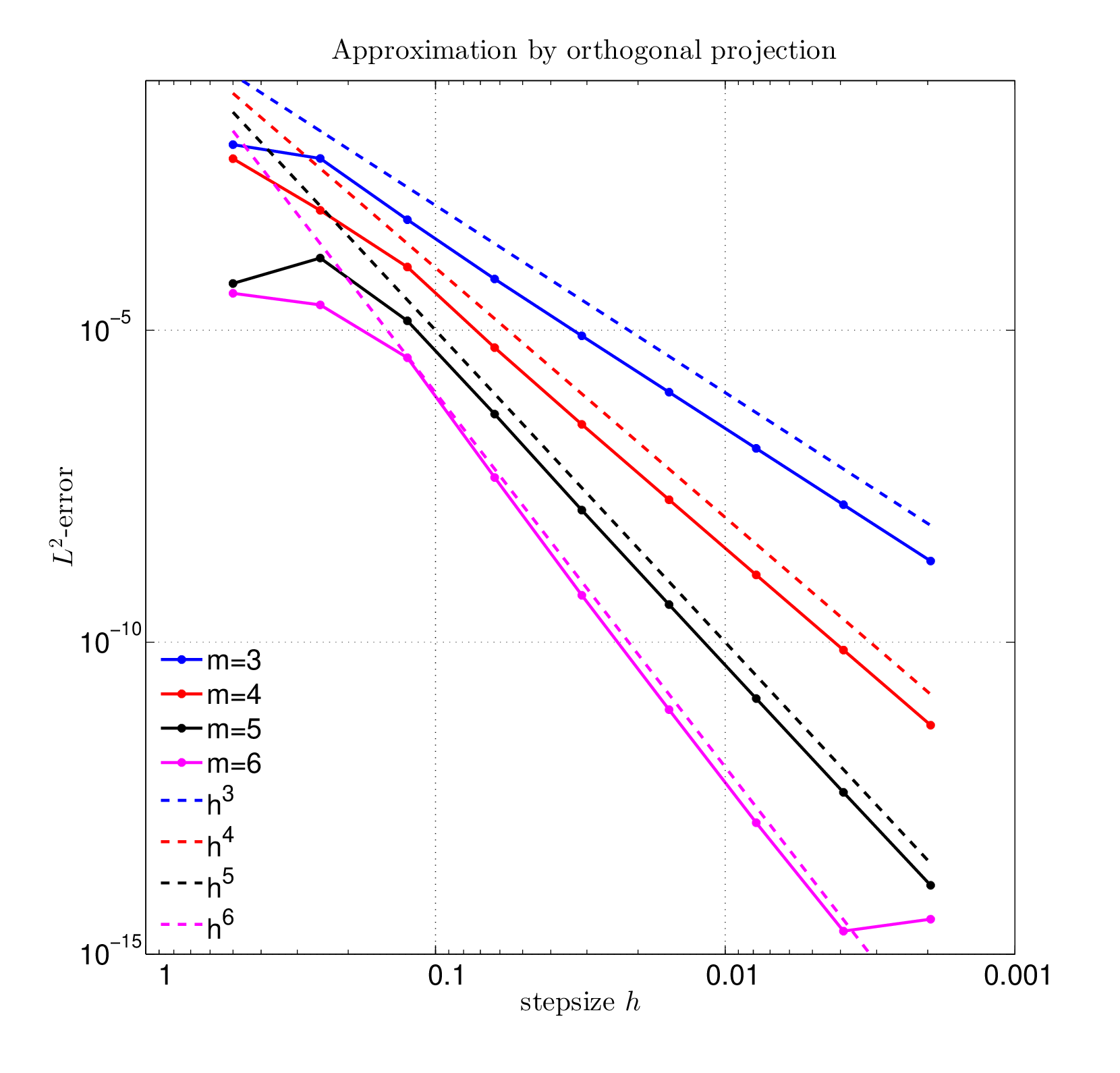}\quad
    \includegraphics[width=0.4\textwidth]{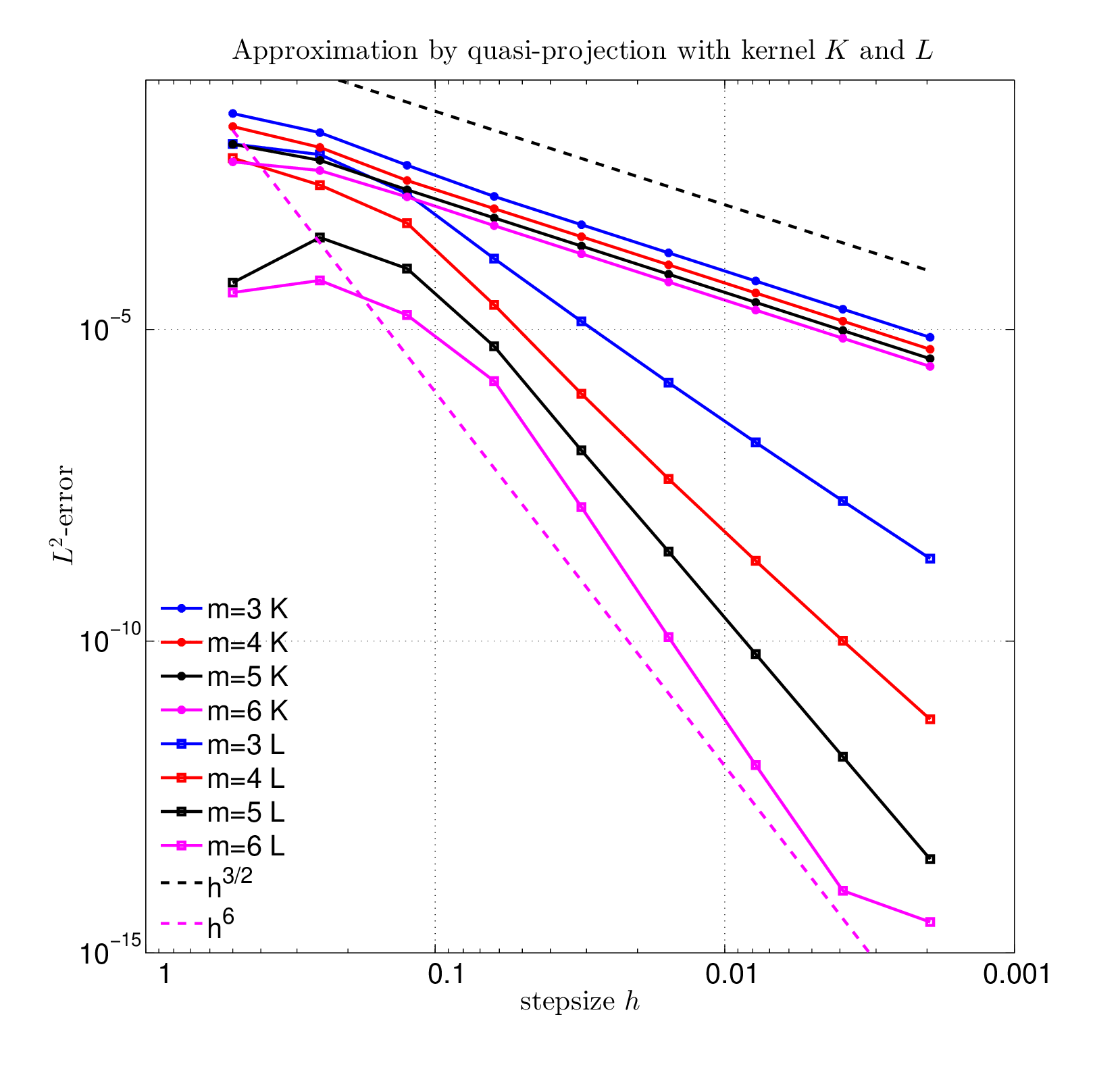}
  \caption{$L^2$-error for approximation of $\hat g$ by splines of order $3\le m\le 6$, with orthogonal projection (left)
  and quasi-projections by means of the kernels $K$ and $L$ (right).
  \label{fig:approxdudn}}
\end{figure}

We choose the knot vector $\Theta_{h,m}^*$ as a refinement of $\Theta_{0,m}^*$ with stepsize $h=\frac{1}{2N}$;
so it is similar to $\Theta_{h,m}$, but has increased multiplicity $m-1$ of the knot $0.5$.
Our numerical results show that
the quasi-projection of $\hat g$ with the kernel $L$ provides the same optimal approximation order as the
orthogonal projection, whereas the kernel $K$ only provides the approximation order $O(1.5)$, as predicted by
the classical Sobolev smoothness of $\hat g$. The results are shown in Figure \ref{fig:approxdudn}.

\newpage

\section{Explicit formulation of enhanced approximate duals}

\noindent
We give explicit formulations of the kernel $L$ for the cases $m=2$ and $m=3$.

\subsection{Case $m=2$}

\noindent
Let $\Theta$ be a knot vector with double knots $\theta_1=\theta_2=a$ and $\theta_ {n+1}= \theta_{n+2}=b$ at both
endpoints.
Since we consider the case $m=2$, all interior knots $\theta_3,\ldots, \theta_{n}$ are simple.

The Gramian of $\Phi_2$ is the tridiagonal matrix
$$ \Gamma = \int_a^b \Phi_2(x)^T \Phi_2(x)\,dx = (\gamma_{j,k}; ~ 1\le j,k\le n) $$
with entries
$$ \gamma_{j,j} = \frac{\theta_{j+2}-\theta_j}{3},\qquad \gamma_{j+1,j}=\gamma_{j,j+1} = \frac{\theta_{j+2}-\theta_{j+1}}{6}.$$
The kernel $K$
in \cite{Chuietal2004} for the reproduction of linear polynomials has the form
\begin{eqnarray*}  K(x,y)&=&   \Phi_2(y)~S~\Phi_2(x)^T   \\
&=& \sum_{k=1}^n \frac{2}{\theta_{k+2}-\theta_k} N_{2,k}(x)N_{2,k}(y) +
\sum_{k=1}^{n-1} \frac{(\theta_{k+2}-\theta_{k+1})^2}{2(\theta_{k+3}-\theta_k)} N_{3,k}'(x)N_{3,k}'(y) ,
\end{eqnarray*}
where we specify the explicit values of the coefficients $u^{(\nu)}_k$ in Eq. \eqref{eq:khq1} as given in
\cite[Eq. (5.13)]{Chuietal2004}. From this expression, we find
the matrix $S$ in \eqref{eq:defS}: it is a tridiagonal matrix with  entries
$$ s_{j,j}= \frac{2}{\theta_{j+2}-\theta_j}+\frac{2}{(\theta_{j+2}-\theta_j)^2}(\alpha_{j-1}+\alpha_j),\qquad  s_{j+1,j}=s_{j,j+1} = \beta_j,$$
where
$$ \alpha_0=\alpha_n=0,\qquad \alpha_j= \frac{(\theta_{j+2}-\theta_{j+1})^2}{\theta_{j+3}-\theta_{j}}\quad\text{for}\quad 1\le j\le n-1$$
and
$$ \beta_j =- \frac{2(\theta_{j+2}-\theta_{j+1})^2}{(\theta_{j+2}-\theta_{j})(\theta_{j+3}-\theta_{j})(\theta_{j+3}-\theta_{j+1})}\quad
 \text{for}\quad  1\le j\le n-1.$$
The product $\Gamma S$ is a matrix with upper and lower bandwidth $2$.
We use symbolic computations with Matlab for the factorization of $I_n-\Gamma S$ and obtain the following result.

\begin{lemma}  We have
$$ I_{n}-\Gamma S = Z D_3 D_2$$
with a lower triangular matrix $Z\in \bR^{n\times (n-2)}$, whose entries are
\begin{eqnarray*} z_{j,j}&=&\frac{(\theta_{j+2}-\theta_{j+1})(\theta_{j+3}-\theta_{j+2})^2}{18(\theta_{j+3}-\theta_{j+1})},\quad
z_{j+1,j}= -\frac{(\theta_{j+2}-\theta_{j+1})(\theta_{j+3}-\theta_{j+2})}{18},\\
z_{j+2,j}&=& \frac{(\theta_{j+2}-\theta_{j+1})^2(\theta_{j+3}-\theta_{j+2})}{18(\theta_{j+3}-\theta_{j+1})}
\end{eqnarray*}
for $1\le j\le n-2$ and whose remaining entries are zero.
\end{lemma}

Our ansatz for the kernel $L$ in \eqref{eq:Lm} is
\begin{equation}\label{eq:ansatz}
  L(x,y)=\Phi_2(y)~S_L~\Phi_2(x)^T \quad\text{with}\quad  S_L = S+  D_2^T\,D_3^T\, U_2\,D_3\,D_2
\end{equation}
and with a banded matrix $U_2$ of size $(n-2)\times (n-2)$.

Let $3\le \ell\le n$ be an index of an interior knot of $\Theta$, for which we require the reproduction of the
truncated power $p_\ell(x):=(\theta_{\ell}-x)_+$. Marsden's identity \eqref{eq:Marsden2b} for truncated powers gives
\begin{equation}\label{eq:cell}
   p_\ell = \Phi_2 c_\ell \quad\text{with} \quad c_\ell = (\theta_{\ell}-\theta_2,\cdots,\theta_{\ell}-\theta_{\ell-1},0,\cdots,0)^T.
\end{equation}
By symbolic computations we obtain
\begin{equation}\label{eq:kappaell}
 c_\ell^T Z = \kappa_\ell e_{\ell-2}^T\quad\text{with}\quad \kappa_\ell=
\frac{(\theta_{\ell}-\theta_{\ell-1})^2(\theta_{\ell+1}-\theta_{\ell})^2}{18(\theta_{\ell+1}-\theta_{\ell-1})}.
\end{equation}
Here, $e_k$ is the canonical unit vector in $\bR^{n-2}$ with $1$ in its $k$-th entry.
This gives the following explicit form for $U_2$ in Proposition \ref{prop:findUm}.

\begin{proposition}
Let $3\le \ell_1 < \cdots < \ell_r\le n$ be indices of interior knots of $\Theta$.
Every solution $U_2$ of the equation $AU_2=B$, where
\begin{equation}\label{eq:U2condition}
      A =\left[ N_{4,k}(\theta_{\ell_j})\right]_{j=1,\ldots,r;~k=1,\ldots,n-2} \in \bR^{r\times (n-2)},\qquad B= \left[\kappa_{\ell_j} e_{\ell_j-2}^T\right]\in \bR^{r\times (n-2)}
\end{equation}
and $\kappa_{\ell_j}$ is given in \eqref{eq:kappaell}, defines a kernel $L$ in \eqref{eq:ansatz}
which reproduces all polynomials of degree $1$ and all truncated powers
$(\theta_{\ell_j}-x)_+$ with $1\le j\le r$.
\end{proposition}

\begin{remark}\label{rem:linearspline}
\rm
The matrix $A_0$ in \eqref{eq:A0} is
$$ A_0 = \left[ N_{4,\ell_s-2}(\theta_{\ell_t})\right]_{s,t=1,\ldots,r} \in \bR^{r\times r}.$$
We restrict the following discussion to the special case where
$$ \ell_j\le \ell_{j+1}-2\quad\text{for all}\quad 1\le j\le r-1.$$
Then $A_0$ is a diagonal matrix, and so is $U_2$. The only non-zero entries of $U_2$ are
$$
(U_2)_{\ell_j-2,\ell_j-2} =  \kappa_{\ell_j} / N_{4,\ell_j-2}(\theta_{\ell_j}),\quad  1\le j\le r .
$$
The explicit form of the kernel $L$ in \eqref{eq:ansatz} is
\[ L(x,y) = K(x,y) +
\sum_{j=1}^r \frac{ \kappa_{\ell_j}}{N_{4,\ell_j-2}(\theta_{\ell_j})} N_{4,\ell_j-2}''(x)N_{4,\ell_j-2}''(y).
\]
\end{remark}

\subsection{Case $m=3$}

\noindent
In a similar way, we can describe the construction of enhanced approximate duals for quadratic B-splines in $\cS_3(\Theta)$.
The knot vector $\Theta$ has triple knots $\theta_1=\theta_2=\theta_3=a$ and $\theta_ {n+1}= \theta_{n+2}=\theta_{n+3}=b$ at both
endpoints.
All interior knots $\theta_4,\ldots, \theta_{n}$ have multiplicity $1$ or $2$ (which would mean $\theta_j=\theta_{j+1}$ for some $j$).
Our ansatz for the new kernel $L$  is
\begin{equation}\label{eq:ansatz2}
  L(x,y)=\Phi_3(y)~S_L~\Phi_3(x)^T \quad\text{with}\quad  S_L = S+  D_3^T\,D_4^T\,D_5^T\, U_3\, D_5\,D_4\,D_3
\end{equation}
with a banded matrix $U_3$ of size $(n-3)\times (n-3)$.
Our subsequent results for $m=3$ are obtained using symbolic software (Matlab). We only present our results for the case of simple knots,
so $\mu_k=1$ for $4\le k\le n$.

\begin{proposition} Let $\Theta$ be a knot vector with simple interior knots $\theta_4,\ldots,\theta_n$ and triple knots at the endpoints.
\begin{itemize}
\item[(a)]  We have
$$ I_{n}-\Gamma S = Z D_5 D_4 D_3$$
with a matrix $Z\in \bR^{n\times (n-3)}$, whose entries $z_{i,j}$ are zero if $j>i+1$; in other words, $Z$ has only one nonzero upper diagonal.
\item[(b)] Let $4\le \ell\le n$ and $p_\ell(x):=(\theta_{\ell}-x)_+^2$.
 Marsden's identity \eqref{eq:Marsden2b} for truncated powers gives
$ p_\ell = \Phi_3 c_\ell $ with
\begin{equation}\label{eq:cell2}
  c_\ell=(c_{\ell,j},~1\le j\le n),\quad  c_{\ell,j} = \begin{cases} (\theta_{\ell}-\theta_{j+1})(\theta_\ell-\theta_{j+2}),&~ 1\le j\le \ell-3,\\
     0,&~\text{otherwise}.  \end{cases}
\end{equation}
We let $(a,b,c,d,e,f)=(\theta_{j+1}-\theta_j;~\ell-3\le j\le \ell+2)$.
Then
\begin{equation}\label{eq:kappaell2a}
 w_\ell := c_\ell^T Z = \mu_\ell e_{\ell-4}^T+ \kappa_\ell e_{\ell-3}^T +\nu_\ell e_{\ell-2}^T,
\end{equation}
where $e_k$ is the canonical unit vector in $\bR^{n-3}$  and the scalar factors are
\begin{equation}\label{eq:kappaell2b}
 \mu_\ell =\frac{(b^2(a+b+c)^2+a (a+b)(b+c)c) d^4 }{600 (a+b+c+d) (b+c+d) (c+d)},
\end{equation}
\begin{equation}\label{eq:kappaell2c}
\nu_\ell =\frac{(e^2 (d+e+f)^2+ f(e+f)(d+e) d) c^4 }{600 (c+d+e+f) (c+d+e)  (c+d)},
\end{equation}
{\overfullrule0pt
\begin{equation}\label{eq:kappaell2d}
\begin{array}{rcl}
  \kappa_\ell &=&\displaystyle \frac{1}{600 (b+c+d) (b+c+d+e) (c+d) (c+d+e)}  \biggl(
 c^2(b + c)(c + d)(d + e)^3(b + c + d)^2 + \\[12pt]
 && c^2 d^3 (b + c) (c + d + e) (b + c + d + e)^2 + c^2 d (b + c) (d + e)^2 (b + c + d)^2 (c + d + e) + \\[8pt]
&& c^2 d^2 (b + c) (c + d) (d + e) (b + c + d + e)^2 + b c^2 d^2 (d + e)^2 (b + c + d) (b + c + d + e) + \\[8pt]
&& b d^2 (b + c)^2 (c + d) (d + e)^2 (c + d + e) +  b^2 c d^2 (d + e)^3 (b + c + d) + \\[4pt]
&& b^2 c d^3 (d + e)^2 (b + c + d + e)\biggr).
\end{array}
\end{equation}
}
$\kappa_\ell=\kappa_\ell(b,c,d,e)$ is symmetric in the sense
$$\kappa_\ell(b,c,d,e)= \kappa_\ell(e,c,d,b)= \kappa_\ell(b,d,c,e) =\kappa_\ell(e,d,c,b).$$
Its numerator is a homogenous polynomial of degree $9$.
Furthermore,
$$ c_\ell^T \, \Gamma\,D_3^T\, D_4^T\, D_5^T = \left( N_{6,j}(\theta_{\ell}); ~1\le j\le n-3\right).$$
\item[(c)] We fix several indices $4\le \ell_1 < \cdots < \ell_r\le n$ of interior knots of $\Theta$. Every solution $U_3$ of the equation $AU_3=B$,
\begin{equation}\label{eq:U2condition2}
      A =2~\left[ N_{6,k}(\theta_{\ell_j})\right]_{j=1,\ldots,r;~k=1,\ldots,n-3} \in \bR^{r\times (n-3)},\qquad B= \left[w_{\ell_j}\right]_{j=1,\ldots,r}\in \bR^{r\times (n-3)}
\end{equation}
with row vectors $w_{\ell_j}$ in \eqref{eq:kappaell2a} defines a kernel $L$ which reproduces all polynomials of degree $2$ and all truncated powers
$(\theta_{\ell_j}-x)_+^2$ with $1\le j\le r$.
\end{itemize}
\end{proposition}

\begin{remark}
{\rm
Under similar conditions as in Remark \ref{rem:linearspline}, the matrix $A_0$ in \eqref{eq:A0} is a diagonal matrix, so taking its inverse
is trivial.
More precisely, this is the case if
\begin{equation}\label{eq:condquadratic}
 \ell_j\le \ell_{j+1}-3\quad\text{for}\quad 1\le j\le r-1.
\end{equation}
Therefore, the use of the enhanced approximate duals will cause almost no computational overhead, if only few truncated powers
should be reproduced in addition to all polynomials of degree $2$.
}
\end{remark}


\end{document}